\documentclass[12pt]{article}

\usepackage{amsmath,amsthm}

\def\v{\vert}
\def\a{\ensuremath{\mathcal A}}

\def\f{\ensuremath{\mathcal F}}
\def\g{\ensuremath{\mathcal G}}

\begin{document}
\newtheorem{lemma}{Lemma}
\newtheorem{theorem}{Theorem}
\newtheorem{prop}{Proposition}

\begin{center}
{\Large
 Pattern avoidance in circular permutations         \\ 
}
\vspace{10mm}
DAVID CALLAN  \\
Department of Statistics  \\
University of Wisconsin-Madison  \\
1210 W. Dayton St   \\
Madison, WI \ 53706-1693  \\
{\bf callan@stat.wisc.edu}  \\
\vspace{5mm}
\today
\end{center}

\vspace{5mm}

In the present context, a \emph{circular permutation} is an 
arrangement of distinct positive 
integers (called letters) clockwise around a circle, 
two arrangements being considered the same 
if they differ only by a rotation. A (linear) permutation is a list of 
distinct positive integers and its \emph{reduced form} is obtained by 
replacing its smallest entry by 1, its next smallest by 2, and so on. 
A \emph{pattern} is any such reduced form. An occurrence of a 
pattern $\tau$ in a cyclic permutation $\pi$ is a sequence of letters 
in clockwise order in $\pi$ (not extending more than one
revolution) whose reduced form is $\tau$. We will 
always represent a circular permutations on $[n]:=\{1,2,\ldots,n\}$ as a linear
permutation that ends in 
$n$. For example, $5\,6\,4\,2\,3\,1\,7$ has just one occurrence of the 
pattern 1234, the letters being $2\,3\,5\,6$ (keep in mind we may wrap around from 
the end to the start).

Much work has been done counting permutations with restrictions on 
the patterns they contain and possibly other restrictions as 
well; see \cite{wilf02} for a survey. There are two 3-letter patterns up to cyclic 
rotation: 123 and 132. Only the identity circular permutation avoids 132 and only 
the reverse identity avoids 123. Here we find simple explicit formulas for the 
number of circular permutations avoiding a single four letter 
pattern. There are 6 such patterns up to cyclic rotation and ``reversing'' 
cuts the number to 3:  1234, 1324 and 1342  are 
representative. 
In each case, an analysis of the structure 
according to the position of $n-1$ (recall $n$ is at the end) yields 
a recurrence relation. For completeness we include a proof of the 
following result.
\newpage 
\begin{prop}[\cite{simion-schmidt}]
	There are $2^{n-1}$ $($linear$)$ permutations on $[n]$ that avoid both $213$ and 
	$231$.
\end{prop}
\begin{proof}
	Consider such a permutation $\pi$ and any letter $a$ in $\pi$. The 
	letters following $a$ in $\pi$ must either all be $>a$ or all be 
	$<a$, else $a$ would be the first letter of a proscribed pattern. So 
	the the first letter of $\pi$ is 1 or $n$ and we have two choices for 
	each successive letter except the last. This observation translates 
	into a bijection from 0-1 sequences of length $n-1$ to 
	permutations on $[n]$ that avoid 213 and 231: given 
	$(b_{i})_{i=1}^{n-1},\ b_{i}=0$ or 1, form $\pi$ left to right as 
	follows. For $i=1$ to $n-1$ in turn, if $b_{i}=1$ (resp. 0), set 
	$\pi_{i}=$ the largest (resp. smallest) letter in $[n]$ that has not 
	already been used in $\pi$. For example, with $n=8,\ (0,1,1,1,0,1,0) \to 
	(1,8,7,6,2,5,3,4)$.
\end{proof}

\begin{theorem}
	The number of $1324$-avoiding circular permutations on $[$n$]$ is 
	the Fibonacci number $F_{2n-3}$
\end{theorem}
\begin{proof}
	Let $u_{n}$ denote the desired number. Then $u_{1}=u_{2}=1,\ 
	u_{3}=2$. For $n\ge 4$, we count by the 
	location $k \in [n-1]$ of $n-1$ in $\pi$. First suppose $ k \le n-2$. 
	so that $n-1$ and $n$ are not adjacent.
	The first $k-1$ letters in $\pi$ must all be $>$ than all the 
	letters between $n-1$ and $n$, else there would be a 1324 pattern of 
	the form $a \ n\!-\!1\ b\ n$. So the first $k$ entries are 
	$\{n-j\}_{j=1}^{k}$. Furthermore, they must occur in increasing 
	order, else there would be a $c \ b\ n\!-\!1\ a$ pattern, 
	rotationally equivalent to a $1324$. 
	We leave the reader to verify that any such circular permutation 
    is in fact 1324-avoiding.
	It is now easy to check that 
	for $k \le n-2$ (including $k=1$), deleting the first $k$ entries and 
	reducing (in this case, reducing means replacing $n$ by $n-k$) is a 
	bijection to 1324-avoiding circular permutations on $[n-k]$, counted 
	by $u_{n-k}$. Now suppose $k=n-1$. Deleting the $n$ gives a bijection 
	to the $n-1$ case, counted by $u_{n-1}$.
	Summing over $k$, we get a recurrence relation: 
	$u_{n}=u_{n-1}+\sum_{k=1}^{n-1}u_{n-k}\ n \ge 3$ with initial 
	condition $u_{1}=u_{2}=1$, whose solution is  
	$u_{n}=F_{2n-3}$.
	
	Here is a short combinatorial proof that $u_{n}:=F_{2n-3}$ satisfies 
	the recurrence. Let $\f_{n}$ denote the set of 0-1 sequences of 
	length $n-1$ that start with 0 (unless they're empty) and do not 
	contain two consecutive 1's, so that $\v \f_{n} \v =F_{n}$. Set 
	$\g_{n}=\f_{2n-3}$ so that $\v \g_{n} \v =u_{n}$. If $w \in \g_{n}$ 
	begins 01\ldots, deleting this 01 is a bijection to $\g_{n-1}$. 
	Otherwise $w$ begins $00\ldots$ and we count by number $k$ of 
	consecutive 10
	pairs immediately following the initial 00. Clearly, $k \in [0,n-3]$ 
	and deleting the initial 00 and all $k$ of these 10's is a bijection to 
	$\g_{n-1-k}$, counted by $u_{n-1-k}$. 
\end{proof}
\newpage 
\begin{theorem}
	There are $2^{n-1}-(n-1)$ $1342$-avoiding circular permutations on $[n]$.
\end{theorem}
\begin{proof}
	Let $u_{n}$ denote the number of 1342-avoiding circular permutations on 
	$[n]$.  Once again $u_{1}=u_{2}=1,\ 
	u_{3}=2$, so let $n \ge 4$. If $n\!-\!1$ is in first 
	position, deleting $n$ (and rotating) gives a bijection to 1342-avoiding 
	circular permutations on 
	$[n-1]$, counted by $u_{n-1}$.
	
	Now suppose $n\!-\!1$ is in position $k$ with $k \in [2,n-2].$ Write $\pi 
	=w_{1}\ n\!-\!1\ \ldots\ n$ with $w_{1}$ the initial subword of 
	length $k-1$. Then $w_{1}$ is nonempty and its letters are 
	increasing, else $n\!-\!1,n$ would be the middle letters of a 1342 
	pattern. Also, the letters in $w_{1}$  form an interval of consecutive 
	integers. To see this, suppose $a<c$ are letters of $w_{1}$ and 
	$b\in[a+1,c-1]$. Then $b$ can't precede $a$ in $\pi$ or fall between 
 	$c$ and $n\!-\!1$ because the letters in $\pi$ increase. Neither can 
 	$c$ fall between $n\!-\!1$ and $n$, else $a\ c\ n\!-\!1\ b$ would be a 
 	$1342$. So $b$ must fall between $a$ and $c$ and $w_{1}$ is a list 
 	of consecutive integers.
	
	It follows that  $w_{1}$ splits the remaining letters in $[n-2]$ into two sets 
	$A,B$ (either or both may be empty, $A$ denotes the set of smaller 
	letters). Write $\pi$ as $w_{1}\ n\!-\!1\ w_{0}\ n$ so that $w_{0}$ is 
	comprised of the letters in $A\cup B$. The letters in $A$ are in increasing 
	order in $\pi$ because if $w_{0}=\ldots\,a_{2}\,\ldots\,a_{1}\,\ldots$ 
	with $a_{1}<a_{2}$, then $a_{1}\,b\ n\!-\!1\ a_{2}$ 
	would be a 1342 for any $b \in 
	w_{1}$. All letters of $C$ precede all 
	letters of $A$ because if $w_{0}=\ldots\,a\,\ldots\,c\,\ldots$ with 
	$a\in A,\,c\in C$, then $a\,c\,n\,b$ would be a 1342 for any $b \in 
	w_{1}$.
	
	So far, then, we have shown that $\pi$ has the form
	\[
    w_{1}\ n\!-\!1\ w_{2}\ w_{3} \ n
	\]
	with $w_{1}$ of length $k-1$ (and so nonempty), all letters in $w_{2}>$ those in $w_{1}$, all letters in 
	$w_{3}$ increasing and $<$ those in $w_{1}$. Now we claim further 
	that the word $w_{2}$ must avoid the patterns 213 and 231. Indeed, if
	$w_{2}=\ldots\,b\,\ldots\,a\,\ldots\,c\,\ldots$ with $a<b<c$, then 
	$a\,c\,n\,b$ is a 1342 in $\pi$. And if 
	$w_{2}=\ldots\,b\,\ldots\,c\,\ldots\,a\,\ldots$, then 
	$x\,b\,c\,a$ is a 1342 for any $x$ in $w_{1}$. We leave the 
	reader to verify that any circular permutation with $n\!-\!1$ in 
	position $k \in [2,n-2]$  
	that meets all these specifications is in fact 1342-avoiding. Let us 
	count these latter permutations. The initial segment $w_{1}$
	is determined by its first letter $j$ which must lie in $[n-k]$. 
	Then $w_{3}$ is fully determined as are the letters of $w_{2}$. 
	Since $w_{2}$ avoids 213 and 231 (Prop.\:1), there are  $2^{\v 
	w_{2}\v-1}=2^{n-k-j-1}$ possibilities for $w_{2}$ unless $j=n-k$ in 
	which case $w_{2}$ is empty and there is just one possibility. 
	Summing over $j$, we have 
	$1+\sum_{j=1}^{n-k-1}2^{n-k-j-1}=2^{n-k-1}$ possibilities in all. 
	Now summing over $k$, we get the defining recurrence $u_{1}=1,\ 
	u_{2}=1,\ u_{n}=u_{n-1}+2^{n-2}-1$ for $n \ge 3$, 
	with solution $u_{n}=2^{n-1}-(n-1)$  valid for all $n \ge 1$.
	
	The preceding analysis could be worked up into a bijection. Briefly, 
	let $\a_{n}$ denote the set of $n$-bit 0-1 sequences that do not 
	contain precisely one 1 so that $\v \a_{n}\v =2^{n}-n$. Given $u 
	\in \a_{n-1}$ build up $\pi=w_{1}\ n\!-\!1\ w_{2}\ w_{3} \ n$ as follows. 
	If $u_{1}=0$, put $n\!-\!1$ in 
	position 1 and proceed recursively. Otherwise, $u_{1}=1$ and put 
	$n\!-\!1$ in the position $k\ge 2$ of the next 1 in $u$ (there will be one!). 
	Then let $j \in [n-k]$ be the 
	position of the next 1 after that in $u_{k+1}\ldots u_{n-1}$ 
	(with $j=n-k$ if there are no more 1's). Use $j$ as the initial letter 
	that determines $w_{1}$ and use 
	$u_{k+j+1}\ldots u_{n-1}$ to construct $w_{2}$ (Prop. 1). 
\end{proof}
\begin{theorem}
	There are $2^{n}+1-2n-\binom{n}{3}$ $1234$-avoiding circular permutations on $[n]$.
\end{theorem}
\begin{proof}
	Again let $u_{n}$ denote the desired number and $k$ the position of 
	$n\!-\!1$. If $k=1$, deleting the $n\!-\!1$ and reducing is a bijection to 
	the $n-1$ case. If $k=2$, write $\pi$ as $a\ n\!-\!1\ w_{1} \ n $ 
	with $a \in [n-2]$. All letters  $>a$ in $w_{1}$ clearly occur in 
	decreasing order, as must all letters $>a$ (else $n\!-\!1$ would 
	terminate a 1234). So the only freedom is in placement of the 
	letters $<a$, giving $\sum_{a=1}^{n-2}\binom{n-3}{a-1}=2^{n-3}$ 
	possibilities with $k=2$. 
	
	Now suppose $k \ge 3$ so that 
	$\pi=w_{1}\ n\!-\!1\ w_{2}\ n$ with $\v w_{1}\v=k-1 \ge 2$. Let $m$ and $M$ 
	denote the minimum and maximum letters in $w_{1}$ and let $A$ 
	consist of the remaining $k-3$ letters in $w_{1}$. Then the letters 
	of $w_{2}$ can be partitioned into three sets: $B=[1,m-1],\ 
	C=[m+1,M-1]\backslash A,\ D=[M+1,n-2]$ of sizes $\v B \v =m-1,\ \v 
	C\v =M-m-1-K, \v D\v =n-M-2$ respectively, where $K=\v A\v=k-3$. All letters $<M$ 
	in $w_{2}$ must occur in decreasing order, else there would be a 
	1234 terminating in $M\ n\!-\!1$; likewise for all letters $>m$, else 
	$muvn$ would be a 1234 for some $u,v\in w_{2}$. This means that 
	unless $M=m+1$, \emph{all} letters in $w_{2}$ must decrease and 
	in case $M=m+1$ (and hence $C$ is empty) the letters of $B$ and $D$ 
	separately must occur in decreasing order in $w_{2}$ giving 
	$\binom{n+m-M-3}{m-1}$ choices to place the letters of $B$ and $D$ 
	in $w_{2}$. Again we leave the reader to verify that if all these 
	conditions are met, the proscribed pattern is avoided.
	
	Using standard binomial coefficient identities, we can now count possibilities. 
	Fix $k \in [3,n-1]$ (so $K=k-3$ is also fixed) and sum over $m$ and $M$. There are 
	$\binom{M-m-1}{K}$ ways to select the set $A$. For $C=\emptyset$ (in 
	which case $M=m+1+K$), the count is 
	$\sum_{m=1}^{n-K-3}\binom{n-K-4}{m-1}=2^{n-K-4}=2^{n-1-k}$. For $C\ne 
	\emptyset$, the count is 
	\begin{eqnarray*}
	\sum_{m=1}^{n-K-4}\sum_{M=m+K+2}^{n-2}\binom{M-m-1}{K} & = &
	\sum_{m=1}^{n-K-4}\left(\binom{n-m-2}{K+1}-1\right) \\
	& = & \binom{n-2}{K+2}-n+K+3 \\
	& = & \binom{n-2}{k-1}-n+k
	\end{eqnarray*}
	Taking the cases $C=\emptyset$ and $C\ne \emptyset$ together, the total 
	count for fixed $k \in [3,n-1]$ is $2^{n-1-k} + 
	\binom{n-2}{k-1}-n+k$, while, as shown above, for $k=2$ it is 
	$2^{n-3}$ and for $k=1$ it is $u_{n-1}$. Summing over $k$ now yields 
	the recurrence relation $u_{1}=u_{2}=1,\ 
	u_{n}=u_{n-1}+2^{n-1}-n-\binom{n-2}{2}$ for $n \ge 3$, with 
	solution $u_{n}=2^{n}+1-2n-\binom{n}{3}$.
\end{proof}

We conclude by noting that the preceding exact enumerations yield growth 
constants ($\lim_{n \to \infty}u_{n}^{1/n}$)  
for the patterns 1324,\ 1342 and 1234 
of $\phi^{2} \approx 2.618,\ 2$ and 2 respectively, where $\phi$ denotes the golden ratio.

\end{document}